\documentclass[11pt]{amsart}

\usepackage{amscd,amsmath,amssymb,fancyhdr,color}
\usepackage[utf8x]{inputenc}
\usepackage{amsfonts}
\usepackage{amsthm}
\usepackage{accents}
\usepackage{graphicx}
\usepackage{float}
\usepackage{subcaption}
\usepackage{verbatim}
\usepackage{indentfirst}
\usepackage{dsfont}
\usepackage{tikz}
\usepackage{enumerate}
\usepackage[all,cmtip]{xy}
\usepackage{faktor}

\usepackage[backref=page]{hyperref}
\renewcommand*{\backref}[1]{}
\renewcommand*{\backrefalt}[4]{%
    \ifcase #1 (Not cited.)%
    \or        (Cited on page~#2.)%
    \else      (Cited on pages~#2.)%
    \fi}

\hypersetup{
     colorlinks   = true,
     citecolor    = magenta
}

\numberwithin{equation}{section}

\def\eqref#1{(\ref{#1})}

\def\1{\sqrt{-1}\:}

\newcommand{\cntrct}                % contraction with a vector field
{\hspace{2pt}\raisebox{1pt}{\text{$\lrcorner$}}\hspace{2pt}}

\newcommand{\codim}{\operatorname{codim}}
\renewcommand{\dim}{\operatorname{dim}}

\newcounter{Mycounter}[section]
\newcounter{lemma}[section]
\newcounter{claim}[section]

\newcounter{sublemma}[section]

\newcounter{corollary}[section]

\newcounter{theorem}[section]

\newcounter{conjecture}[section]

\newcounter{proposition}[section]

\newcounter{definition}[section]

\newcounter{example}[section]

\newcounter{remark}[section]

\newcounter{problem}[section]

\newcounter{question}[section]

\makeatletter

\@addtoreset{equation}{section}

\@addtoreset{footnote}{section}

\makeatother

\usetikzlibrary{arrows,chains,matrix,positioning,scopes}

\makeatletter
\tikzset{join/.code=\tikzset{after node path={%
			\ifx\tikzchainprevious\pgfutil@empty\else(\tikzchainprevious)%
			edge[every join]#1(\tikzchaincurrent)\fi}}}
\makeatother

\tikzset{>=stealth',every on chain/.append style={join},
	every join/.style={->}}

\DeclareMathOperator{\proj}{proj}

\begin{document}

\title[$q$-pseudoconvex and $q$-holomorphically convex domains]
 {$q$-pseudoconvex and $q$-holomorphically convex domains}

%----------Author 1
\author{George-Ionu\c t Ioni\c t\u a}

\address{%
\textsc{Department of Mathematics and Computer Science\\ 
"Politehnica" University Bucharest}\\ 
313 Splaiul Independen\c tei, Bucharest 060042, Romania}

\email{georgeionutionita@gmail.com}

%----------Author 2
\author{Ovidiu Preda}

\address{%
\textsc{Institute of Mathematics of the Romanian Academy}\\
P.O. Box 1-764, Bucharest 014700, Romania}

\email{ovidiu.preda@imar.ro}

\thanks{Ovidiu Preda was
supported by a grant of Ministry of Research and Innovation, CNCS - UEFISCDI, project number
PN-III-P1-1.1-PD-2016-0182, within PNCDI III}

%----------classification, keywords, date

\begin{abstract}
In this article we prove a global result in the spirit of Basener's theorem
regarding the relation between $q$-pseudoconvexity and $q$-holomorphic convexity:
we prove that any smoothly bounded strictly $q$-pseudoconvex open subset $\Omega\subset\mathbb{C}^n$ is 
$(q+1)$-holomorphically convex; moreover, assuming that $\Omega$ verifies an additional assumption, we prove that it is $q$-holomorphically convex.
We also prove that any open subset of $\mathbb{C}^n$ is 
$n$-holomorphically convex.
\end{abstract}

%%% ----------------------------------------------------------------------
\maketitle
%%% ----------------------------------------------------------------------
%\tableofcontents
\section{Introduction}

Basener \cite{BAS} proved that the solutions of the equation $\overline{\partial}f\wedge(\partial
\overline{\partial}f)^{q-1}=0$, where $f$ is a smooth function on a complex manifold and $q\geq 1$ is an  integer, define a notion of convexity which is, locally, related to $q$-pseudoconvexity 
in the same way that holomorphic convexity is related to pseudoconvexity. He proved \cite[Theorem 3]{BAS} 
the following result: 
\textit{
Let $\Omega\subset\mathbb{C}^n$ be a smoothly bounded, bounded, open subset. If $\Omega$ is $q$-holomorphically 
convex, then it is $q$-pseudoconvex. Also, if $\Omega$ is strictly $q$-pseudoconvex, then it is locally 
$q$-holomorphically convex.
}

In this article we prove a global $q$-holomorphic convexity result, similar 
to \cite[Theorem 3, (ii)]{BAS}, but 
with either weaker convexity conclusion or stronger hypothesis. Our first result is the following:

\begin{theorem}\label{TH1}
Let $\Omega$ be a smoothly bounded, strictly $q$-pseudoconvex domain in $\mathbb{C}^n$. Then, $\Omega$ is 
$(q+1)$-holomorphically convex. Moreover, if $\Omega$ is bounded and 
for every $p\in\partial\Omega$, there exists a closed complex submanifold 
$X$ of $\mathbb{C}^n$ of dimension $n-q+1$, which is a holomorphic complete intersection, 
such that $p\in X$, $X\cap\Omega$ is strictly pseudoconvex, and $X$ intersects $\partial\Omega$ transversally, then $\Omega$ is $q$-holomorphically convex.
\end{theorem}

Also, Greene and Wu \cite{GW} proved that every non-compact complex manifold of dimension 
$n$ is $n$-complete. In particular, any open subset of $\mathbb{C}^n$ is $n$-complete. Here we prove
a similar result, for $n$-holomorphic convexity:

\begin{theorem}\label{TH2}
Any open subset of $\mathbb{C}^n$ is $n$-holomorphically convex.
\end{theorem}

\section{Preliminaries}

In this section we collect the definitions and results needed for \ref{TH1} and \ref{TH2}.
For the coherence of the notations, we use the convention introduced by Andreotti and Grauert \cite{AG} 
for defining $q$-convexity and we adapt all the other definitions according to it, even if 
Basener \cite{BAS} has used a different convention.

Let $\Omega$ be an open subset of $\mathbb{C}^{n}$ and $f:\Omega \rightarrow \mathbb{R}$ a smooth function (for simplicity, ``smooth" will stand for $\mathcal{C}^\infty$). If $z_{0} \in \Omega$, then the Levi form of $f$ at $z_{0}$, denoted $L(f,z_{0})$, is the quadratic form determined by
$(\partial^{2}f/\partial z_{i}\partial \bar{z}_{j}(z_{0}))_{i,j}$.
 A function $f$ is called $q$-convex if its Levi form has at least $n-q+1$ positive ($>0$) eigenvalues at any point of $\Omega$.

We observe that the definition of $q$-convex function (which in~\cite{AG} is called strictly $q$-pseudoconvex) coincides with the one of a strictly $(q-1)$-convex function given by Basener~\cite{BAS}.

Consider now $\Omega$ to be a smoothly bounded domain in $\mathbb{C}^n$. Let $x\in\partial\Omega$ and let $U$ be an open 
neighborhood of $x$ in $\mathbb{C}^n$. A smooth real-valued function $\phi$ on $U$ is called a defining function on $U$ 
for $\Omega$ if $\Omega\cap U=\{\phi<0\}$ and $d\phi\not=0$ on $\partial\Omega\cap U$. If the restriction of the Levi form 
$L(\phi,z)$ to the tangent space $T_p(\partial\Omega)$ has at least $n-q$ positive (respectively, nonnegative) eigenvalues, 
then $\Omega$ is said to be strictly $q$-pseudoconvex (respectively, $q$-pseudoconvex) at $p$.
As usual, $\Omega$ is said to be strictly $q$-pseudoconvex (respectively, $q$-pseudoconvex) if it is strictly $q$-pseudoconvex 
(respectively, $q$-pseudoconvex) 
at each boundary point. 
1-pseudoconvexity will be called just pseudoconvexity.

\begin{definition}
Let $\Omega$ be a complex manifold. Define 
$$\mathcal{O}_q(\Omega)=\{f\in\mathcal{C}^\infty(\Omega): \overline{\partial}f\wedge(\partial\overline{\partial}f)^{q-1}=0\}.$$
If $f\in\mathcal{O}_q(\Omega)$, we will say that $f$ is $q$-holomorphic.
\end{definition}
It is easy to see that 1-holomorphic functions are exactly the functions which are holomorphic in the classical sense.

\begin{definition}
 Let $\Omega$ be a complex manifold and let $K$ be a compact subset of $\Omega$.
The set 
$$\widehat{K}_{\mathcal{O}_q(\Omega)}=\{z\in\Omega: |f(z)|\leq \max_K|f| \ \mbox{for all}\  
f\in \mathcal{O}_q(\Omega) \}$$
is called the $q$-holomorphically convex hull of $K$. We will say that $\Omega$ is $q$-holomorphically 
convex if for each compact set $K\subset \Omega$, the set $\widehat{K}_{\mathcal{O}_q(\Omega)}$ 
is again compact.
\end{definition}

\begin{definition}
 Let $D$ be a domain in $\mathbb{C}^n$ with smooth boundary. We denote by $A^\omega(D)$ the space 
of functions which are holomorphic on a neighborhood of $\overline{D}$. We say that a point $p\in\partial D$ 
is a (global) peak point relative to $D$ for $A^\omega(D)$ if there exists a function $f\in A^\omega(D)$ 
such that $f(p)=1$ and $|f|<1$ on $\overline{D}\setminus \{p\}$. We call $f$ a peak function. We say that 
$p\in\partial D$ is a local peak point for $A^\omega(D)$ if $p$ is a peak point for $A^\omega(D\cap V)$ 
for some neighborhood $V$ of $p$.
\end{definition}

The following theorem by Rossi \cite[Theorem 4.4]{ROS}, also mentioned in Noell's survey 
\cite[Remark 2.4 (2)]{NOE}, gives a sufficient condition for the existence of peak functions:

\begin{theorem}\label{NOE_THM}
 If $D$ is an open subset of a Stein manifold and $\overline{D}$ has a basis of strictly pseudoconvex domains, then every local 
peak point for $A^\omega(D)$ is a global peak point for $A^\omega(D)$. 
In particular, if $D$ is relatively compact and strictly 
pseudoconvex at every boundary point, then $\overline{D}$ has a basis of strictly pseudoconvex domains, 
so every boundary point is a global peak point for $A^\omega(D)$.
\end{theorem}

The next theorem, commonly known as Narasimhan's lemma, states that every strictly pseudoconvex boundary 
point of a smoothly bounded domain has holomorphic local coordinates in which it is strictly convex.

\begin{theorem}\label{NAR_LEM}
 Let $D\subset \mathbb{C}^n$ be a smoothly bounded domain which is strictly pseudoconvex at a point 
$p\in\partial D$. 
Then, there exists a biholomorphic map $F$ defined in a neighborhood $U$ of $p$ such that 
$F(D\cap U)$ is strictly convex in $F(U)$.
\end{theorem}
The following theorem by Docquier and Grauert (\cite{DOC-GRA}, \cite[p.257, Theorem 8]{GUN-ROS}) states the existence of tubular neighborhoods for complex submanifolds. 
Although their result is stronger than what we need and we will use only a part of it, we give here 
the original statement:
\begin{theorem}\label{TH_DG}
Let $S$ be a Stein submanifold of a complex manifold $X$. Denote by $N_{S/X}$ the normal bundle of $S$ in $X$. Then, there exists an open Stein neighborhood $U$ of $S$ in $X$, biholomorphic to an open neighborhood $\Omega$ of the zero section in $N_{S/X}$, and a homotopy of holomorphic maps $\iota_t:U\rightarrow U$ $(t\in[0,1])$ such that $\iota_0$ is the identity map on $U$, $\iota_t|_S$ is the identity map on $S$ for all $t\in[0,1]$, and $\iota_1(U)=S$.
\end{theorem}

\begin{definition}
 A closed complex submanifold $Y$ of codimension $d$ in a complex manifold $X$ is a holomorphic complete 
intersection if there exist $d$ holomorphic functions $f_1,\ldots,f_d\in\mathcal{O}(X)$ such that 
$Y=\{x\in X:f_1(x)=\ldots=f_d(x)=0\}$, and the differentials $df_j(x)$ $(1\leq j\leq d)$ are $\mathbb{C}$-linearly 
independent at each point $x\in Y$.
\end{definition}

It is easy to prove that these differentials induce a trivialization of the normal bundle $N_{Y/X}={T_X}_{|Y}/T_Y$, leading 
to the following lemma:

\begin{lemma}\label{NORMB}
 If $Y$ is a closed complex submanifold of a complex manifold $X$, such that $Y$ is a holomorphic 
complete intersection, then the normal bundle $N_{Y/X}$ of $Y$ in $X$ is trivial.
\end{lemma}

\section{The main results}
\subsection{Proof of \ref{TH1}}

\

\textit{Strategy of the proof.} Since $\mathcal{O}_q(\Omega)$ is not closed under addition, there seems to be no way to patch together $q$-holomorphic 
functions, so a global theorem cannot be derived from the local result proved by Basener. The plan for our proof is the 
following: for a given point on the boundary of $\Omega$, we construct a convenient closed submanifold $X\subset\mathbb{C}^n$ 
passing through that point and whose 
intersection with $\Omega$ is strictly pseudoconvex; then we use the existence of peak holomorphic functions on this 
intersection, which 
we extend differentiably to $\Omega$ by making them go to zero in the normal directions given by the closed submanifold $X$;
in this way, we obtain ''almost peak $(q+1)$-holomorphic functions`` which are sufficient to show the $(q+1)$-holomorphic convexity 
of $\Omega$. With the additional assumption on $\Omega$ stated in \ref{TH1}, the previously mentioned closed submanifold $X$ 
is already given by the hypothesis and has higher dimension, and this allows us to reduce $(q+1)$ to $q$.

It is easy to observe that the proof we give for the first assertion in \ref{TH1} can be shorter. However, we prefer to write this proof in a form which, with the additional assumption, also solves the second assertion.

\begin{proof}

For the first part of the proof, we know that $\Omega$ is a smoothly bounded, 
strictly $q$-pseudoconvex domain in $\mathbb{C}^n$, and we want to show that $\Omega$ is 
$(q+1)$-holomorphically convex.

We prove that for each $p\in\partial\Omega$ and each neighborhood $V_p$ of $p$, there exists a function 
$f\in\mathcal{C}^\infty(\overline{\Omega})$ with $f_{|\Omega}\in\mathcal{O}_{q+1}(\Omega)$,
satisfying $f(p)=1$, $|f|<1$ on $\overline{\Omega}\setminus V_p$. The existence of these 
$(q+1)$-holomorphic functions is sufficient to prove that $\Omega$ is 
$(q+1)$-holomorphically convex. If $K\subset \Omega$ is a compact subset, these functions 
show that the closure of 
$\widehat{K}_{\mathcal{O}_{q+1}(\Omega)}$ in $\mathbb{C}^n$ does not contain points of $\partial \Omega$. 
Also, it is easy to see that $\widehat{K}_{\mathcal{O}_{q+1}(\Omega)}$ is closed in $\Omega$ and 
bounded in $\mathbb{C}^n$. Therefore, $\widehat{K}_{\mathcal{O}_{q+1}(\Omega)}$ is compact, yielding that 
$\Omega$ is 
$(q+1)$-holomorphically convex.

We begin now the construction of these $(q+1)$-holomorphic functions. Fix $p\in\partial\Omega$ and 
$V_p$ a 
neighborhood of $p$.
Since $\Omega$ is strictly $q$-pseudoconvex at $p$, there exists a neighborhood $U$ of $p$ and a 
defining function $\phi:U\rightarrow \mathbb{R}$ on $U$ for $\Omega$ such that the Levi form 
$L(\phi,p)$ restricted to the tangent space $T_{p}(\partial\Omega)$ has at least $n-q$ positive eigenvalues.
Denote by $L$ a complex linear space through $p$, of dimension $n-q+1$, 
which contains the direction normal to the boundary of $\Omega$ in $p$, namely $\nabla \phi(p)$, and 
other $n-q$ directions in the tangent space $T_p(\partial\Omega)$ which are directions 
of positivity for the Levi form $L(\phi,p)$ (we know there are 
at least $n-q$). Then, $\Omega\cap L$ is an open 
set, with smooth boundary near $p$, which has a strict pseudoconvex boundary point at $p$.

Next, by \ref{NAR_LEM}, there exists a neighborhood $U$ of $p$ in $L$ and a 
biholomorphism $F:U\rightarrow F(U)\subset \mathbb{C}^{n-q+1}$ such that $F(\Omega\cap L\cap U)$ 
is strictly convex in $F(U)$. We may assume, in addition, that $F(U)$ is a ball in $\mathbb{C}^{n-q+1}$. 
Now, we take a linear complex subspace $M$ in $\mathbb{C}^{n-q+1}$ of dimension $n-q$, containing $p$, such 
that $M\cap F(\Omega\cap L\cap U)$ is nonempty and relatively compact in $F(U)$. 
Since $F(\Omega\cap L\cap U)$ is strictly convex in $F(U)$, it follows that 
$M\cap F(\Omega\cap L\cap U)$ is also strictly convex, hence strictly pseudoconvex. 
We denote by $\proj_L:\mathbb{C}^n\rightarrow L$ the projection on the linear subspace $L$.
Then, $F^{-1}(M)$ is a Stein manifold contained in the open subset 
$D:=\proj_L^{-1}(U)\subset\mathbb{C}^n$ 
and $F^{-1}(M)\cap (\Omega\cap L\cap U)$ is a strictly pseudoconvex, relatively compact, 
open subset of $F^{-1}(M)$, with $p$ one of its boundary points in $F^{-1}(M)$. 

$M$ is linear, so $M$ is a holomorphic complete intersection in $F(U)$. Since $F$ is a biholomorphism, 
$F^{-1}(M)$ is a holomorphic complete intersection in $U$; hence, knowing that $\codim_U M=1$, we obtain that there exists 
$f_1\in\mathcal{O}(U)$ such that $F^{-1}(M)=\{z\in U:f_1(z)=0\}$ and $df_1(z)\not =0$ for all $z\in F^{-1}(M)$.
We can extend $f_1$ holomorphically to $\widetilde{f}_1:D=\proj_L^{-1}(U)\rightarrow \mathbb{C}$ by making it 
constant on each fiber.
Also, $L$ is a complex linear subspace of $\mathbb{C}^n$ of codimension $q-1$, so it is a holomorphic complete intersection: 
$L=\{z\in\mathbb{C}^n: l_1(z)=\ldots=l_{q-1}(z)=0\}$, where $l_1=\ldots=l_{q-1}$ are linear functions and $dl_1,\ldots,dl_{q-1}$ 
are linearly independent on $L$. Since $l_1=\ldots=l_{q-1}$ are zero on $L$ and $\widetilde{f}_1$ is constant on the fibers 
of the projection on $L$, we obtain that 
$d\widetilde{f}_1,dl_1,\ldots,dl_{q-1}$ are linearly independent on $F^{-1}(M)$.
Consequently, 
$F^{-1}(M)=\{z\in D:\widetilde{f}_1(z)=l_1(z)=\ldots=l_{q-1}(z)=0 \}$ is a holomorphic complete intersection in $D$.

By \ref{NOE_THM}, there exists an open, relatively compact neighborhood $W$ of $F^{-1}(M)\cap (\Omega\cap L\cap U)$ in $F^{-1}(M)$ and a holomorphic function $h:W\rightarrow \mathbb{C}$ such that $h(p)=1$ and $|h|<1$ on $\overline{F^{-1}(M)\cap (\Omega\cap L\cap U)}\setminus\{p\}$. 

By \ref{TH_DG}, $F^{-1}(M)$ has a neighborhood in $D$ which is biholomorphic 
to an open neighborhood of the zero section in the normal bundle $N_{F^{-1}(M)/D}$. Since 
$F^{-1}(M)$ is a holomorphic complete intersection in $D$, \ref{NORMB} implies that the bundle 
$N_{F^{-1}(M)/D}$ is trivial. Hence, we have $N_{F^{-1}(M)/D}\simeq F^{-1}(M)\times \mathbb{C}^{q}$.
Since $W$ is relatively compact in $F^{-1}(M)$, the biholomorphism given by \ref{TH_DG} can be restricted 
to $G:Q\subset\mathbb{C}^n\rightarrow \overline{W}\times B(0,r)$. 
We may assume that we have chosen both $W\subset F^{-1}(M)$ and $r>0$ small 
enough such that $G^{-1}(\partial W\times B(0,r))\cap \Omega=\emptyset$ and 
$G^{-1}(\{w\in W:|h(w)|\geq 1\})\times B(0,r))\subset V_p$.

We consider now a smooth decreasing function $g:\mathbb{R}\rightarrow [0,1]$ such that 
$g(t)=1$ for $t\leq 0$ and $g(t)=0$ for $t\geq r$. 

Finally, we define $u:W\times B(0,r)$, $u(w,b)=h(w)g(\|b\|)$. Since $u$ is holomorphic in the $n-q$ 
variables corresponding to $W$, $u$ is $(q+1)$-holomorphic. Consider $f:\Omega\rightarrow \mathbb{C}$, 
$f(z)=u(G(z))$ if $z\in Q$ and $f=0$ on $\Omega\setminus Q$. Then, $f$ is a smooth function which is 
$(q+1)$-holomorphic, $f(p)=1$ and $|f|<1$ on $\overline{\Omega}\setminus V_p$. 
Hence, as mentioned in the beginning of the proof, we can conclude that $\Omega$ is 
$(q+1)$-holomorphically convex. 

For the ``moreover'' part of the conclusion, the proof is exactly the same, with the closed complex submanifold $X$ which has the 
properties mentioned in the theorem's statement, 
instead of $F^{-1}(M)$. Since $\dim X=\dim F^{-1}(M) +1$, we obtain that $\Omega$ is $q$-holomorphically 
convex.
\end{proof}

Since for any $q$, the class of $q$-holomorphically convex open subsets is closed under finite intersections, \ref{TH1} provides new examples of $q$-holomorphically convex subsets:
\begin{corollary}
 Any finite intersection of strictly $q$-pseudoconvex domains in $\mathbb{C}^n$ is $(q+1)$-holomorphically convex.
\end{corollary}

\subsection{Proof of \ref{TH2}}

\

Basener \cite[Example 5, p.205]{BAS} proves that for every $\lambda=(\lambda_1,\ldots,\lambda_n)\in\mathbb{C}^n\setminus\{0\}$, the function $f_\lambda:\mathbb{C}^n\setminus\{0\}\rightarrow \mathbb{C}$ defined by 
$$f_\lambda(z_1,\ldots,z_n)=\left(\sum_{i=1}^{n}\lambda_i\overline{z}_i\right)\left(\sum_{i=1}^{n}|z_i|^2\right)^{-1}$$
is $n$-holomorphic and has an isolated nonremovable singularity at the origin. We use this type of functions to prove \ref{TH2}.

\begin{proof}
Consider $K\subset \Omega$ a compact subset. 
It is easy to see that $\widehat{K}_{\mathcal{O}_n(\Omega)}$ is 
closed in $\Omega$ and bounded in $\mathbb{C}^n$.
To prove that $\widehat{K}_{\mathcal{O}_n(\Omega)}$ is compact, it is enough to show that 
the boundary of $\widehat{K}_{\mathcal{O}_n(\Omega)}$ in $\mathbb{C}^n$ does not contain points of 
$\partial\Omega$.

It is sufficient to show that for every $p\in\partial\Omega$ and every $r>0$ such that 
$B(p,r)\cap K=\emptyset$, we have
$B(p,\frac{r}{\sqrt{n}})\cap \widehat{K}_{\mathcal{O}_n(\Omega)}=\emptyset$. 
For this, we fix $p\in\partial\Omega$ and $r>0$ such that $B(p,r)\cap K=\emptyset$.
Then, consider $z=(z_1,\ldots,z_n)\in B(p,\frac{r}{\sqrt{n}})$ and take 
$\lambda=(\lambda_1,\ldots,\lambda_n)\in 
\mathbb{C}^n$ such that $|\lambda_1|=\ldots=|\lambda_n|=1$ and 
$\lambda_i \overline{(z_i-p_i)}=|z_i-p_i|$ for any 
$1\leq i\leq n$.  Then, for any $w\in\partial B(p,r)$, we have: 
$$|f_\lambda(z-p)|=\frac{\sum|z_i-p_i|}{\|z-p\|^2}\geq \frac{1}{\|z-p\|}>
\frac{\sqrt{n}}{\|w-p\|}\geq \frac{\sum|w_i-p_i|}{\|w-p\|^2}\geq$$
$$ 
\geq \left|\frac{\sum \lambda_i \overline{(w_i-p_i)}}{\|w-p\|^2}\right|=|f_\lambda(w-p)|
$$
Also, the function $(0,\infty)\ni t\mapsto |f_\lambda(tx)|=\frac{1}{t}|f_\lambda(x)|$ 
is strictly decreasing. Hence, $\displaystyle \max_{x\in K}|f_\lambda(x-p)|\leq\max_{x\in\partial B(p,r)}|f_\lambda(x-p)|<
|f_\lambda(z-p)|$, 
so $z\not\in \widehat{K}_{\mathcal{O}_n(\Omega)}$. Since $z\in B(p,\frac{r}{\sqrt{n}})$ was chosen arbitrarily, we get that 
$B(p,\frac{r}{\sqrt{n}})\cap \widehat{K}_{\mathcal{O}_n(\Omega)}=\emptyset$, which ends our proof.

\end{proof}

We remark that using global holomorphic embeddings in the Euclidean space,
 \ref{TH2} remains true, with the same proof, for Stein manifolds instead of $\mathbb{C}^n$.

\section{Remarks}

The additional assumptions in \ref{TH1} lead us to the following problem:

\begin{problem}\label{PB1}
 Let $\Omega\subset \mathbb{C}^n$ be an open, strictly $q$-pseudoconvex subset, where $1< q\leq n$, and $p\in\partial \Omega$. 
Does there exist a closed submanifold $X$ of $\mathbb{C}^n$, of dimension $n-q+1$, which is a holomorphic complete intersection, 
such that $p\in X$, $X$ intersects 
$\partial\Omega$ transversally, and $X\cap\Omega$ is strictly pseudoconvex?
\end{problem}

An affirmative solution to this problem would lead to the following sequence of implications, for bounded domains with 
smooth boundary: strictly $q$-pseudoconvex $\Rightarrow$ $q$-holomorphically convex $\Rightarrow$ $q$-pseudoconvex. 
However, even if the additional condition from \ref{TH1} is not true in general, 
for many given domains it can be easily checked that it is verified. 
A weaker version of \ref{PB1} can be stated as follows:

\begin{problem}\label{PB2}
 Let $\Omega\subset \mathbb{C}^n$ be an open, strictly $q$-pseudoconvex subset, where $1< q\leq n$, and $p\in\partial \Omega$. 
Do there exist $f_1,\ldots,f_{q-1}\in\mathcal{O}(\mathbb{C}^n)$, such that $p\in X=\{f_1=\ldots=f_{q-1}=0\}$, and 
$X\cap\Omega$ is Stein?
\end{problem}

% ------------------------------------------------------------------------

% ------------------------------------------------------------------------
\end{document}